\documentclass[11pt]{article}
\usepackage{amsmath}
\usepackage{amssymb}
\usepackage{amsbsy}
\usepackage{amsfonts}
\usepackage{mathtools}
\usepackage{bbm}
\usepackage{comment}
\newtheorem{theorem}{Theorem}[section]%
\newtheorem{lemma}[theorem]{Lemma}%
\newtheorem{cor}[theorem]{Corollary}%
\newenvironment{pf}{\medskip\noindent{Proof:}
\hspace{-.5cm} \enspace}{\hfill \qed \newline \smallskip}

\def\f{\noindent}

\newcommand{\qed}{\mbox{\raisebox{0.7ex}{\fbox{}}} \vspace{4truemm}}

\begin{document}

\begin{center}
{\Large{\textbf{On Jacobian group and complexity of  $I$-graph $I(n,k,l)$ through Chebyshev polynomials}}}
\end{center}

\vskip 5mm {\small
\begin{center}{\textbf{I.~A.~Mednykh}}\footnote{{\small\em Sobolev Institute of Mathematics,
Novosibirsk State University, 630090, Novosibirsk, Russia,\\ E-mail: ilyamednykh@mail.ru}}
\end{center}}

\title{\vspace{-1.2cm}
On Jacobian group and complexity of the $I$-graph $I(n,k,l)$ through Chebyshev polynomials
\thanks{Supported by}}

\begin{abstract}
We consider a family of $I$-graphs $I(n,k,l),$ which is a generalization of the class of
generalized Petersen graphs. In the present paper, we provide a new method for counting
Jacobian group of the $I$-graph $I(n,k,l).$ We show that the minimum number of generators
of $Jac(I(n,k,l))$ is at least two and at most $2k+2l-1.$ Also, we obtain a closed formula
for the number of spanning trees of $I(n,k,l)$ in terms of Chebyshev polynomials. We
investigate some arithmetical properties of this number and its asymptotic behaviour.
\bigskip

\f\textbf{Key Words:} Spanning tree, Jacobian group, $I$-graph, Petersen graph, Chebyshev polynomial 
\end{abstract}

\section{Introduction}
The notion of the Jacobian group of a graph, which is also known as the Picard group,
the critical group, and the dollar or sandpile group, was independently introduced by
many authors (\cite{CoriRoss}, \cite{BakerNorine}, \cite{Biggs}, \cite{BachHarpNagnib}).
This notion arises as a discrete version of the Jacobian in the classical theory of Riemann
surfaces. It also admits a natural interpretation in various areas of physics, coding theory,
and financial mathematics. The Jacobian group is an important algebraic invariant of a finite
graph. In particular, its order coincides with the number of spanning trees of the graph, which
is known for some simplest graphs, such as the wheel, fan, prism, ladder, and M\"obius ladder
\cite{BoePro}, grids \cite{NP04}, lattices \cite{SW00}, prism and anti-prism \cite{SWZ16}.
At the same time, the structure of the Jacobian is known only in particular cases \cite{CoriRoss},
\cite{Biggs}, \cite{Lor}, \cite{YaoChinPen}, \cite{ChenHou}, \cite{MedZind} and \cite{MedMed2}.
We mention that the number of spanning trees for circulant graphs is expressed is terms of
the Chebyshev polynomials; it was found in \cite{ZhangYongGol}, \cite{ZhangYongGolin}, and
\cite{XiebinLinZhang}. We show that similar results are also true for the $I$-graph $I(n,k,l).$

The generalized Petersen graph $GP(n,k)$ has vertex set and edge set given by
\begin{eqnarray*}
V(P(n,k)) &=& \{ u_i, v_i \ | \  i=1,2, \ldots, n \} \\
E(P(n,k)) &=& \{ u_{i}u_{i+1}, \ u_iv_i, \ v_i v_{i+k}\ | \ i=1,2, \ldots, n \},
\end{eqnarray*}
where the subscripts are expressed as integers modulo $n$. The classical Petersen graph is $P(5,2)$. The
family of generalized Petersen graphs is a subset of so-called $I$-graphs. The $I$-graph $I(n,k,l)$ is a
graph of the following structure
\begin{eqnarray*}
V(I(n,k,l)) &=& \{ u_i, v_i \ | \  i=1,2, \ldots, n \} \\
E(P(n,k,l)) &=& \{ u_{i}u_{i+l}, \ u_iv_i, \ v_i v_{i+k}\ | \ i=1,2, \ldots, n \}.
\end{eqnarray*}
where all subscripts are given modulo $n.$

Since $I(n,k,l)=I(n,l,k)$ we will usually assume that $k\leq l.$ In this paper we will
deal with $3$-valent graphs only. This means that in the case of even $n$ and $l=n/2$
the graph under consideration has multiple edges. If $GCD(n,k,l)=m>1,$ then $I(n,k,l)$
is a union of $m$ copies of the graph $I(n/m,k/m,l/m).$ If $m=1$ and $GCD(k,l)=d,$ then
the graphs $I(n,k,l)$ and $I(n,k/d,l/d)$ are isomorphic \cite{BobPisZit}, \cite{HorPisZit},
\cite{PetkZakr}. So, in what follows, we assume $k$ and $l$ to be relatively prime. In
the case of $l=1$ it easy to see that the graph $I(n,k,1)$ coincides with the generalized
Petersen graph $GP(n,k).$ The number of spanning trees and the structure of Jacobian group
for the generalized Petersen graph were investigated in \cite{KwonMedMed}. The spectrum
of the $I$-graph was found in \cite{OlivVina}. Even though the number of spanning
trees of a given graph can be computed through eigenvalues of its Laplacian matrix, it is
not easy to find the number of spanning trees for $I(n,k,l)$ using them. In this paper,
we obtained a closed formula for the number of spanning trees for $I(n,k,l),$ investigate
some its arithmetical properties of this number and provide its asymptotic behavior. Also,
we suggest an effective way for calculating Jacobian of $I(n,k,l)$ and find sharp upper and
lower bounds for the rank of $Jac(I(n,k,l)).$

\section{Basic definitions and preliminary facts}

Consider a connected finite graph $G,$ allowed to have multiple edges but without
loops. We endow each edge of $G$ with the two possible directions. Since $G$ has
no loops, this operation is well defined. Let $O=O(G)$ be the set of directed edges
of $G.$ Given $e\in O(G),$ we denote its initial and terminal vertices by $s(e)$
and $t(e),$ respectively. Recall that a closed directed path in $G$ is a sequence
of directed edges $e_i\in O(G),\, i=1,\ldots,n$ such that $t(e_i)=s(e_{i+1})$ for
$i=1,\,\ldots,n-1$ and $t(e_n)=s(e_1).$

Following \cite{BakerNorine} and \cite{BachHarpNagnib}, \textit{the Jacobian group},
or simply \textit{Jacobian } $Jac(G)$ of a graph $G$ is defined as the (maximal) abelian
group generated by flows $\omega(e),e\in O(G),$ obeying the following two Kirchhoff laws:

$K_1:$ the flow through each vertex of $G$ vanishes, that is $\sum\limits_{e\in O,t(e)=x}\omega(e)=0
\textrm{ for all }x\in V(G);$

$K_2:$ the flow along each closed directed path $W$ in $G$ vanishes, that is $\sum\limits_{e\in W}\omega(e)=0.$

\noindent Equivalent definitions of the group $Jac(G)$ can be found in papers \cite{CoriRoss},
\cite{BakerNorine}, \cite{Biggs}, \cite{BachHarpNagnib}, \cite{Lor}, \cite{Dhar}, \cite{KotaniSunada}.

We denote the vertex and edge set of $G$ by $V(G)$ and $E(G),$ respectively. Given $u,v\in V(G),$
we set $a_{uv}$ to be equal to the number of edges between vertices $u$ and $v.$ The matrix
$A =A(G) = \{a_{uv}\}_{u, v\in V (G)},$ called \textit{the adjacency matrix} of the graph $G.$
The degree $d(v)$ of a vertex $v \in V(G)$ is defined by $d(v) = \sum_u a_{uv}.$ Let $D = D(G)$
be the diagonal matrix indexed by the elements of $V(G)$ with $d_{vv} = d(v).$
Matrix $L = L(G) = D(G) - A(G)$ is called \textit{the Laplacian matrix}, or simply \textit{Laplacian},
of the graph $G.$

Recall \cite{Lor} the following useful relation between the structure of the Laplacian matrix and
the Jacobian of a graph $G.$ Consider the Laplacian $L(G)$ as a homomorphism $\mathbb{Z}^{|V|}\to\mathbb{Z}^{|V|},$
where $|V|=|V(G)|$ is the number of vertices in $G.$ The cokernel
$\textrm{coker}\,(L(G))=\mathbb{Z}^{|V|}/\textrm{im}\,(L(G))$ --- is an abelian group. Let
$$\textrm{coker}\,(L(G))\cong\mathbb{Z}_{d_{1}}\oplus\mathbb{Z}_{d_{2}}\oplus\cdots\oplus\mathbb{Z}_{d_{|V|}}$$
be its Smith normal form satisfying the conditions $d_i\big{|}d_{i+1},\,(1\le i\le|V|).$ If the graph
is connected, then the groups $\mathbb{Z}_{d_{1}},\mathbb{Z}_{d_{2}},\ldots,\mathbb{Z}_{d_{|V|-1}}$ ---
are finite, and $\mathbb{Z}_{d_{|V|}}=\mathbb{Z}.$ In this case,
$$Jac(G)\cong\mathbb{Z}_{t_{1}}\oplus\mathbb{Z}_{t_{2}}\oplus\cdots\oplus\mathbb{Z}_{d_{|V|-1}}$$
is the Jacobian of the graph $G.$ In other words, $Jac(G)$ is isomorphic to the torsion subgroup of
the cokernel $\textrm{coker}\,(L(G)).$

Let $M$ be an integer $n\times n$ matrix, then we can interpret $M$ as a homomorphism from $\mathbb{Z}^n$
to $\mathbb{Z}^n.$ In this interpretation $M$ has a kernel $\textrm{ker}\,M,$ an image $\textrm{im}\,M,$ and
a cokernel $\textrm{coker}\,M=\mathbb{Z}^n/\textrm{im}\,M.$ We emphasize that $\textrm{coker}\,M$ of the matrix
$M$ coincides with its Smith normal form.

In what follows, by $I_n$ we denote the identity matrix of order $n.$

We call an $n\times n$ matrix {\it circulant,} and denote it by $circ(a_0, a_1,\ldots,a_{n-1})$ if it is
of the form
$$circ(a_0, a_1,\ldots, a_{n-1})=
\left(\begin{array}{ccccc}
a_0 & a_1 & a_2 & \ldots & a_{n-1} \\
a_{n-1} & a_0 & a_1 & \ldots & a_{n-2} \\
  & \vdots &   & \ddots & \vdots \\
a_1 & a_2 & a_3 & \ldots & a_0\\
\end{array}\right).$$

Recall \cite{PJDav} that the eigenvalues of matrix $C=circ(a_0, a_1,\ldots, a_{n-1})$
are given by the following simple formulas $\lambda_j=p(\varepsilon^j_n),\,j=0,1,\ldots,n-1$
where $p(x)=a_0+a_1 x+\ldots+a_{n-1}x^{n-1}$ and $\varepsilon_n$ is the order $n$ primitive
root of the unity. Moreover, the circulant matrix $C=p(T),$ where $T=circ(0,1,0,\ldots,0)$
is the matrix representation of the shift operator
$T:(x_0, x_1,\ldots,x_{n-2}, x_{n-1})\rightarrow(x_1, x_2,\ldots, x_{n-1},x_0).$

By (\cite{GS11}, lemma 2.1) the $2n\times 2n$ adjacency matrix of the $I$-graph
$I(n,k,l)$ has the following block form
$$A(I(n,k,l)) =
\left(\begin{array}{cc}
C_n^k & I_n\\
I_n & C_n^l\\
\end{array}\right),$$
where $C_n^k$ is the $n\times n$ circulant matrix of the form
$C_n^k = circ(\underbrace{0,\ldots,0}_{k\textrm{ times}},1,0,
\ldots,0,1,\underbrace{0,\ldots,0}_{k-1\textrm{ times}}).$

Denote by $L=L(I(n,k,l))$ the Laplacian of $I(n,k,l).$ Since the graph $I(n,k,l)$ is three-valent, we have
$$L=3I_{2n}-A(I(n,k,l))=\left(\begin{array}{cc}
3I_n-C_n^k & -I_n\\
-I_n & 3I_n- C_n^l\\
\end{array}\right).$$

\smallskip

\section{Cokernels of linear operators}

Let $P(z)$ be a bimonic integer Laurent polynomial. That is $P(z)=z^p+a_1z^{p+1}+\ldots+a_{s-1}z^{p+s-1}+z^{p+s}$
for some integers $p,a_1,a_2,\ldots,a_{s-1}$ and some positive integer $s.$ Introduce the following
companion matrix $\mathcal{A}$ for the polynomial $P(z):$
$\mathcal{A}=\left(\begin{array}{c}\begin{array}{c|c} 0 &
I_{s-1}\end{array}\\\hline -1,-a_1,\ldots,-a_{s-1} \\
\end{array}\right),$ where $I_{s-1}$ is the identity $(s-1)\times(s-1)$ matrix.
We will use the following properties of $\mathcal{A}.$ Note that $\det\mathcal{A}=(-1)^{s}.$
Hence $\mathcal{A}$ is invertible and inverse matrix $\mathcal{A}^{-1}$ is also
integer matrix. The characteristic polynomial of $\mathcal{A}$ coincides with $z^{-p}P(z).$

Let $\mathbb{A}=\langle\alpha_j,\,j\in\mathbb{Z}\rangle$ be a free Abelian group freely
generated by elements $\alpha_j,\,j\in\mathbb{Z}.$ Each element of $\mathbb{A}$ is a
linear combination $\sum\limits_{j}c_j \alpha_j$ with integer coefficients $c_j.$

Define the shift operator $T:\mathbb{A}\rightarrow\mathbb{A}$ as a $\mathbb{Z}$-linear
operator acting on generators of $\mathbb{A}$ by the rule $T:\alpha_j\rightarrow\alpha_{j+1},\,j\in\mathbb{Z}.$
Then $T$ is an endomorphism of $\mathbb{A}.$ Let $P(z)$ be a arbitrary Laurent polynomial
with integer coefficients, then $A=P(T)$ is also an endomorphism of $\mathbb{A}.$ Since
$A$ is a linear combination of powers of $T,$ the action of $A$ on generators $\alpha_j$
can be given by the infinite set of linear transformations
$A:\alpha_j\to\sum\limits_{i}a_{ij}\alpha_i,\, j\in\mathbb{Z}.$ Here all sums under
consideration are finite. We set $\beta_i=\sum\limits_{j}a_{ij}\alpha_j.$ Then $\textrm{im}\,A$
is a subgroup of $\mathbb{A}$ generated by $\beta_i,\,i\in\mathbb{Z}.$ Hence,
$\textrm{coker}\,A=\mathbb{A}/\textrm{im}\,A$ is an abstract Abelian group
$\langle x_i, i\in\mathbb{Z}|\,\sum\limits_{i}a_{ij}x_i=0,\,j\in\mathbb{Z}\rangle$
generated by $x_i,\,i\in\mathbb{Z}$ with the set of defining relations
$\sum\limits_{i}a_{ij}x_i=0,\,j\in\mathbb{Z}.$ Here $x_j$ are images of $\alpha_j$
under the canonical homomorphism $\mathbb{A}\rightarrow\mathbb{A}/\textrm{im}\,A.$
Since $T$ and $A=P(T)$ commute, subgroup $\textrm{im}\,A$ is invariant under the
action of $T.$ Hence, the actions of $T$ and $A$ are well defined on the factor
group $\mathbb{A}/\textrm{im}\,A$ and are given by $T:x_j\rightarrow x_{j+1}$ and
$A:x_j\rightarrow\sum\limits_{i}a_{ij}x_i$ respectively.

This allows to present the group $\mathbb{A}/\textrm{im}\,A$ as follows
$\langle x_i,\,i\in\mathbb{Z}|\,P(T)x_j=0,\,j\in\mathbb{Z}\rangle.$ In a similar
way, given a set $P_1(z),P_2(z),\ldots,P_s(z)$ of Laurent polynomials with integer
coefficients, one can define the group
$\langle x_i,\,i\in\mathbb{Z}|\,P_1(T)x_j=0,P_2(T)x_j=0,\ldots,P_s(T)x_j=0,\,j\in\mathbb{Z}\rangle.$

We will use the following lemma.

\bigskip
\begin{lemma}\label{lemma1} Let $T:\mathbb{A}\rightarrow\mathbb{A}$ be the shift operator.
Consider endomorphisms $A$ and $B$ of the group $\mathbb{A}$ given by the formulas
$A=P(T), B=Q(T),$ where $P(z)$ and $Q(z)$ are Laurent polynomials with integer coefficients.
Then $B:\mathbb{A}\rightarrow\mathbb{A}$ induces an endomorphism $B|_{\textrm{coker}\,A}$
of the group $\textrm{coker}\,A=\mathbb{A}/\textrm{im}\,A$ defined by
$B|_{\textrm{coker}\,A}(\alpha + Im A) = B(\alpha)+Im A,\,\alpha \in \mathbb{A}.$ Furthermore
$$\langle x_i,\,i\in\mathbb{Z}|\,A(T)x_j=0,\,B(T)x_j=0,\,j\in\mathbb{Z}\rangle
\cong\textrm{coker}\,A/\textrm{im}(B|_{\textrm{coker}\,A})\cong\textrm{coker}\,(B|_{\textrm{coker}\,A}).$$\end{lemma}

\begin{pf} The images $\textrm{im}\,A$ and $\textrm{im}\,B$ are subgroups in $\mathbb{A}.$
Denote by $\langle\textrm{im}\,A,\textrm{im}\,B\rangle$ the subgroup generated by elements
of $\textrm{im}\,A$ and $\textrm{im}\,B.$ Since $P(z)$ and $Q(z)$ are Laurent polynomials,
the operators $A=P(T)$ and $B=Q(T)$ do commute. Hence, subgroup $\textrm{im}\,A$ is invariant
under endomorphism $B.$ Indeed for any $y=Ax\in\textrm{im}\,A$, we have $By=B(Ax)=A(Bx)\in\textrm{im}\,A.$
This means that $B:\mathbb{A}\rightarrow\mathbb{A}$ induces an endomorphism of the group
$\textrm{coker}\,A=\mathbb{A}/\textrm{im}\,A.$ We denote this endomorphism by $B|_{\textrm{coker}\,A}.$
We note that the Abelian group $\langle x_i,\,i\in\mathbb{Z}|\,A(T)x_j=0,\,B(T)x_j=0,\,j\in\mathbb{Z}\rangle$
is naturally isomorphic to $\mathbb{A}/\langle\textrm{im}\,A,\textrm{im}\,B\rangle.$ So we have
$$\mathbb{A}/\langle\textrm{im}\,A,\textrm{im}\,B\rangle
\cong(\mathbb{A}/\textrm{im}\,A)/\textrm{im}\,(B|_{\textrm{coker}\,A})
\cong\textrm{coker}\,A/\textrm{im}(B|_{\textrm{coker}\,A})\cong\textrm{coker}\,(B|_{\textrm{coker}\,A}).$$

The lemma is proved.
\end{pf}

\section{Jacobian group for the $I$-graph $I(n,k,l)$}

In this section we prove one of the main results of the paper. We start with the following theorem.

\begin{theorem}\label{theorem0} Let $L=L(I(n,k,l))$ be the Laplacian of the $I$-graph
$I(n,k,l).$ Then $$\textrm{coker}\,L\cong \textrm{coker}(\mathcal{A}^{n}-I),$$ where $\mathcal{A}$ is
$2(k+l)\times 2(k+l)$ companion matrix for the Laurent polynomial $$(3-z^{k}-z^{-k})(3-z^{l}-z^{-l})-1.$$
\end{theorem}

\textbf{Proof.} Let $L$ be the Laplacian matrix of the graph $I(n,k,l).$ Then, as it was mentioned above,
$L$ is a $2n \times 2n$ matrix of the form
$$L=\left(\begin{array}{cc}3I_n-C_n^k & -I_n\\-I_n & 3I_n-C_n^l\\\end{array}\right),$$
where $C_n^k = circ(\underbrace{0,\ldots,0}_{k\textrm{ times}},1,0,\ldots,0,1,
\underbrace{0,\ldots,0}_{k-1\textrm{ times}}).$

Consider $L$ as a $\mathbb{Z}-$linear operator $L:\mathbb{Z}^{2n}\rightarrow\mathbb{Z}^{2n}.$
In this case, $\textrm{coker}(L)$ is an abstract Abelian group generated by elements
$x_1,x_2,\ldots,x_{n},y_1,y_2,\ldots,y_{n}$ satisfying the system of linear equations
$3 x_j-x_{j-k}-x_{j+k}-y_j=0, 3 y_j-y_{j-l}-y_{j+l}-x_j=0$ for any $j=1,\ldots,n.$
Here the indices are considered modulo $n.$ By the property mentioned in Section $2,$
the Jacobian of the graph $I(n,k,l)$ is isomorphic to the finite part of cokernel of
the operator $L$.

To study the structure of $\textrm{coker}(L)$ we extend the list of generators to the
two bi-infinite sequences of elements $(x_j)_{j\in\mathbb{Z}}=(\ldots,x_{-1},x_0,x_1,x_2,\ldots)$
and $(y_j)_{j\in\mathbb{Z}}=(\ldots,y_{-1},y_0,y_1,y_2,\ldots)$ setting $x_{j+m n}=x_{j}$
and $y_{j+m n}=y_{j}$ for any $m\in\mathbb{Z}.$ Then we have the following representation
for cokernel of $L:$

\begin{eqnarray*}\textrm{coker}(L)&=&\langle x_i,y_i,i\in\mathbb{Z}\large| 3 x_j-x_{j+k}-x_{j-k}-y_j=0,\\
&{}&3 y_j-y_{j+l}-y_{j-l}-x_j=0, x_{j+n}=x_{j}, y_{j+n}=y_{j}, j\in\mathbb{Z}\rangle.\end{eqnarray*}

Let $T$ be the shift operator defined by the rule
$T:x_j\rightarrow x_{j+1},\,y_j\rightarrow y_{j+1},j\in\mathbb{Z}.$
Consider the operator $P(T)$ defined by $P(T)=(3-T^{k}-T^{-k})(3-T^l-T^{-l})-1.$ We use the operator
notation from section $3$ to represent the cokernel of $L.$ Then we have
\begin{eqnarray*}
\textrm{coker}(L)
&=&\langle x_i,y_i,i\in\mathbb{Z}\Large|(3-T^k-T^{-k})x_j=y_j,
(3-T^l-T^{-l})y_j=x_j,T^{n}x_j=x_j,T^{n}y_j=y_j,j\in\mathbb{Z}\rangle\\
&=&\langle x_i,i\in\mathbb{Z}\Large|(3-T^l-T^{-l})(3-T^k-T^{-k})x_j=x_j, T^{n}x_j=x_j,j\in\mathbb{Z}\rangle\\
&=&\langle x_i,i\in\mathbb{Z}\Large|((3-T^k-T^{-k})(3-T^l-T^{-l})-1)x_j=0, (T^{n}-1)x_j=0,j\in\mathbb{Z}\rangle \\
&=&\langle x_i,i\in\mathbb{Z}\Large| P(T)x_j=0, (T^{n}-1)x_j=0,j\in\mathbb{Z}\rangle.\\
\end{eqnarray*}

To finish the proof, we apply  Lemma~\ref{lemma1} to the operators
$A=P(T)$ and $B=Q(T)=T^{n}-1.$

Since the Laurent polynomial $P(z)=(3-z^k-z^{-k})(3-z^l-z^{-l})-1$ is bimonic it can be
represented in the form $P(z)=z^{-k-l}+a_1z^{-k-l+1}+\ldots+a_{2k+2l-1}z^{k+l-1}+z^{k+l},$
where $a_1,a_2,\ldots,a_{2k+2l-1}$ are integers. Then the companion matrix $\mathcal{A}$ is
$\left(\begin{array}{c}
\begin{array}{c|c}0 &I_{2k+2l-1}\end{array}\\ \hline
-1,-a_1,\ldots,-a_{2k+2l-1}
\end{array}\right).$ It is easy to see that $\det \mathcal{A}=1$ and its inverse $\mathcal{A}^{-1}$
is also integer matrix.

For convenience we set $s=2k+2l$ to be the size of matrix $\mathcal{A}.$

Note that for any $j\in\mathbb{Z}$ the relations $P(T)x_j=0$ can be rewritten as
$x_{j+s}=-x_{j}-a_1x_{j+1}-\cdots-a_{s-1}x_{j+s-1}.$ Let $\textbf{x}_j=(x_{j+1},x_{j+2},\ldots,x_{j+s})^t$
be $s$-tuple of generators $x_{j+1},x_{j+2},\ldots,x_{j+s}.$ Then the relation $P(T)x_j=0$
is equivalent to $\textbf{ x}_{j}=\mathcal{A}\,\textbf{x}_{j-1}.$ Hence, we have
$\textbf{x}_{1}=\mathcal{A}\,\textbf{x}_0$ and $\textbf{x}_{-1}=\mathcal{A}^{-1}\,\textbf{x}_0,$
where $\textbf{x}_0=(x_{1},x_{2},\ldots,x_{s})^t.$ So, $\textbf{x}_{j} = \mathcal{A}^j\,\textbf{x}_0$
for any $j\in\mathbb{Z}.$ Conversely, the latter implies $\textbf{x}_{j} = \mathcal{A}\,\textbf{x}_{j-1}$
and, as a consequence, $P(T)x_j=0$ for all $j\in\mathbb{Z}.$

Let $\mathbb{A}=\langle\alpha_j,\,j\in\mathbb{Z}\rangle$ be the Abelian group
freely generated by elements $\alpha_j,\,j\in\mathbb{Z}.$ As in Lemma~\ref{lemma1},
we consider operator $A=P(T)$ as an endomorphism of the group $\mathbb{A}.$ Then
$\textrm{coker}\,A=\mathbb{A}/\textrm{im}\,A$ as an abstract Abelian group has the following
representation $\langle\bar{x}_i,i\in\mathbb{Z}\Large| P(T)\bar{x}_j=0,j\in\mathbb{Z}\rangle.$
Here $\bar{x}_j$ are images of generators $\alpha_j$ under the canonical homomorphism
$\mathbb{A}\rightarrow\mathbb{A}/\textrm{im}\,A.$

Our present aim is to show that $\textrm{coker}\,A\cong\mathbb{Z}^s.$ Then we describe
the action of the endomorphism $B|_{\textrm{coker}\,A}$ on the $\textrm{coker}\,A.$
Now we have the following representation of $\textrm{coker}\,A.$
\begin{eqnarray*}
&&\textrm{coker}\,A=\langle\bar{x}_i,i\in\mathbb{Z}\Large| P(T)\bar{x}_j=0,j\in\mathbb{Z}\rangle=\\
&&=\langle \bar{x}_j,j\in\mathbb{Z}\Large|
\bar{x}_{\ell}+a_1\bar{x}_{\ell+1}+\ldots+a_{s-1}\bar{x}_{\ell+s-1}+\bar{x}_{\ell+s}=0,\ell\in\mathbb{Z}\rangle\\
&&=\langle \bar{x}_j,j\in\mathbb{Z}\Large|(\bar{x}_{\ell+1},\bar{x}_{\ell+2},\ldots,\bar{x}_{\ell+s})^t=
\mathcal{A}(\bar{x}_{\ell},\bar{x}_{\ell+1},\ldots,\bar{x}_{\ell+s-1})^t,\ell\in\mathbb{Z}\rangle\\
&&=\langle \bar{x}_j,j\in\mathbb{Z}\Large|(\bar{x}_{\ell+1},\bar{x}_{\ell+2},\ldots,\bar{x}_{\ell+s})^t=
\mathcal{A}^{\ell}(\bar{x}_{1},\bar{x}_{2},\ldots,\bar{x}_{s})^t,\ell\in\mathbb{Z}\rangle\\
&&=\langle \bar{x}_1,\bar{x}_2,\ldots,\bar{x}_s{\Large|}\emptyset \rangle\cong\mathbb{Z}^s.
\end{eqnarray*}

Since the operators $A=P(T)$ and $T$ commute, the action
$T|_{\textrm{coker}\,A}:\bar{x}_j\to \bar{x}_{j+1},\,j\in\mathbb{Z}$ on the $\textrm{coker}\,A$ is
well defined. Now we describe the action of $T|_{\textrm{coker}\,A}$ on the set of generators
$\bar{x}_1,\bar{x}_2,\ldots,\bar{x}_s.$ For any $i=1, \ldots, s-1$, we have
$T|_{\textrm{coker}}(\bar{x}_i)=\bar{x}_{i+1}$ and $T|_{\textrm{coker}\,A}(\bar{x}_s)=\bar{x}_{s+1}=
-\bar{x}_1-a_1\bar{x}_2-\ldots-a_{s-2}\bar{x}_{s-1}-a_{s-1}\bar{x}_s$. Hence, the action of
$T|_{\textrm{coker}\,A}$ on the $\textrm{coker}\,A$ is given by the matrix $\mathcal{A}.$
Considering $\mathcal{A}$ as an endomorphism of the $\textrm{coker}\,A,$ we can write
$T|_{\textrm{coker}\,A}=\mathcal{A}.$ Finally, $B|_{\textrm{coker}\,A}=Q(T|_{\textrm{coker}\,A})=Q(\mathcal{A}).$
Applying Lemma~\ref{lemma1}, we finish the proof of the theorem. $\hfill \qed$

\medskip
\begin{cor}\label{corollary1} The Jacobian group $\textrm{Jac}(I(n,k,l))$ of the $I$-graph
$I(n,k,l)$ is isomorphic to the torsion subgroup of $\textrm{coker}(\mathcal{A}^{n}-I),$ where
$\mathcal{A}$ is the companion matrix for the Laurent polynomial $(3-z^{k}-z^{-k})(3-z^{l}-z^{-l})-1.$
\end{cor}

The Corollary~\ref{corollary1} gives a simple way to find Jacobian group $\textrm{Jac}(I(n,k,l))$
for small values of $k$ and sufficiently large numbers $n.$ The numerical results are given in the
Tables $2$ and $3.$

\section{Counting the number of spanning trees for the $I$-graph $I(n,k,l)$}

\begin{theorem}\label{maintheorem}
The number of spanning trees of the $I$-graph $I(n,k,l)$ is given by the formula
$$\tau_{k,l}(n)=(-1)^{(n-1)(k+l)}n\prod_{s=1}^{k+l-1}\frac{T_n(w_s)-1}{w_s-1},$$
where $w_s, s=1,2,\ldots,k+l-1$ are roots of the order $k+l-1$ algebraic equation
$$\frac{(3-2T_k(w))(3-2T_l(w))-1}{w-1}=0,$$ and $T_j(w)$ is the Chebyshev polynomial
of the first kind.
\end{theorem}

\textbf{Proof.} By the celebrated Kirchhoff theorem, the number of spanning trees
$\tau_{k,l}(n)$ is equal to the product of nonzero eigenvalues of the Laplacian of
a graph $I(n,k,l)$ divided by the number of its vertices $2n.$ To investigate the
spectrum of Laplacian matrix we note that matrix $C^k_n=T^{k}+T^{-k},$ where $T=circ(0,1,\ldots,0)$
is the $n \times n$ shift operator. The latter equality easily follows from the identity $T^n=I_n.$ Hence,
$$ L=\left(\begin{array}{cc}
3I_n-T^{k}-T^{-k}  & -I_n\\
-I_n & 3I_n-T^{l}-T^{-l} \\
\end{array}\right).$$

The eigenvalues of circulant matrix $T$ are $\varepsilon_n^j,$ where
$\varepsilon_n=e^\frac{2\pi i}{n}.$ Since all eigenvalues of $T$ are
distinct, the matrix $T$ is conjugate to the diagonal matrix
$\mathbb{T}=diag(1,\varepsilon_n,\ldots,\varepsilon_n^{n-1})$, where
diagonal entries of
$diag(1,\varepsilon_n,\ldots,\varepsilon_n^{n-1})$ are
$1,\varepsilon_n,\ldots,\varepsilon_n^{n-1}$. To find spectrum of
$L,$ without loss of generality, one can assume that $T=\mathbb{T}.$
Then the $n \times n$ blocks of $L$ are diagonal matrices. This
essentially simplifies the problem of finding eigenvalues of $L.$
Indeed, let $\lambda$ be an eigenvalue of $L$ and
$(x,y)=(x_1,\ldots,x_n,y_1,\ldots,y_n)$ be the corresponding
eigenvector. Then we have the following system of equations

$$\left\{\begin{array}{cc}
(3I_n-T^{k}-T^{-k})x-y & =\lambda x\\
-x+(3I_n-T^{l}-T^{-l})y & = \lambda y \\
\end{array}.\right.$$

From here we conclude that $y=(3I_n-T^{k}-T^{-k})x-\lambda x =(3-\lambda-T^{k}-T^{-k})-1)x.$
Substituting $y$ in the second equation, we have $((3-\lambda-T^{l}-T^{-l})(3-\lambda-T^{k}-T^{-k})-1)x=0$.

Recall the matrices under consideration are diagonal and the $(j+1,j+1)$-th entry of $T$ is equal to
$\varepsilon_n^{j}.$ Therefore, we have
$((3-\lambda-\varepsilon_n^{jk}-\varepsilon_n^{-jk})(3-\lambda-\varepsilon_n^{jl}-\varepsilon_n^{-jl})-1)x_{j+1}=0$
and $y_{j+1}=(3-\lambda-\varepsilon_n^{jl}-\varepsilon_n^{-jl})x_{j+1}.$

So, for any $j=0,\ldots, n-1$ the matrix $L$ has two eigenvalues, say $\lambda_{1,j}$ and $\lambda_{2,j}$
satisfying the quadratic equation
$(3-\lambda-\varepsilon_n^{jk}-\varepsilon_n^{-jk})(3-\lambda-\varepsilon_n^{jl}-\varepsilon_n^{-jl})-1=0.$
The corresponding eigenvectors are $(x,y),$ where $x=\textbf{e}_{j+1}=(0,\ldots,\underbrace{1}_{(j+1)-th},\ldots,0)$
and $y=(3-\lambda-T^{k}-T^{-k})\textbf{e}_{j+1}$. In particular, if $j=0$ for $\lambda_{1,0}, \lambda_{2,0}$
we have $(1-\lambda)(1-\lambda)-1=\lambda(\lambda-2)=0.$ That is, $\lambda_{1,0}=0$ and $\lambda_{2,0}=2.$
Since $\lambda_{1,j}$ and $\lambda_{2,j}$ are roots of the same quadratic equation, we obtain
$\lambda_{1,j}\lambda_{2,j}=P(\varepsilon_n^j),$ where $P(z)=(3-z^{k}-z^{-k})(3-z^{l}-z^{-l})-1.$

Now we have $$\tau_{k,l}(n)=\frac{1}{2n}\lambda_{2,0}\prod\limits_{j=1}^{n-1}\lambda_{1,j}\lambda_{2,j}=
\frac{1}{n}\prod\limits_{j=1}^{n-1}\lambda_{1,j}\lambda_{2,j}=\frac{1}{n}\prod\limits_{j=1}^{n-1}P(\varepsilon_n^j).$$

To continue we need the following lemma.

\begin{lemma}\label{lemma2}
The following identity holds $$(3-z^{k}-z^{-k})(3-z^{l}-z^{l})-1=(3-2T_k(w))(3-2T_l(w))-1,$$
where $T_k(w)$ is the Chebyshev polynomial of the first kind and $w=\frac{1}{2}(z+z^{-1}).$
Moreover, if $k$ and $l$ are coprime then all roots of the Laurent polynomial
$(3-z^{k}-z^{-k})(3-z^{l}-z^{-l})-1$ counted with multiplicities are
$1,\,1,\,z_1,\,1/z_1,\ldots,\,z_{k+l-1},\,1/z_{k+l-1},$ where we have $|z_s|\neq1,\,s=1,2,\ldots,k+l-1.$
So, the right-hand polynomial has the roots $1,\,w_1,\ldots,\,w_{k+l-1},$ where $w_s\neq1$
for all $s=1,\,2,\ldots,k+l-1.$
\end{lemma}

\textbf{Proof.}
Let us substitute $z=e^{i\,\varphi}.$ It is easy to see that $w=\frac{1}{2}(z+z^{-1})=\cos\varphi,$
so we have $T_k(w)=\cos(k\arccos w)=\cos(k\varphi)$. Then the first statement of the lemma is
equivalent to the following trigonometric identity
$$(3-2\cos(k\varphi))(3-2\cos(l\varphi))-1=(3-2T_k(w))(3-2T_l(w))-1.$$

To prove the second statement of the lemma we suppose that the Laurent polynomial
$P(z)=(3-z^{k}-z^{-k})(3-z^{l}-z^{-l})-1$ has a root $z_0$ such that $|z_0|=1.$
Then $z_0=e^{i\,\varphi_0},\,\varphi_0\in\mathbb{R}.$ Now we have
$(3-2\cos(k\varphi_0))(3-2\cos(l\varphi_0))-1=0.$ Since $3-2\cos(k\varphi_0)\geq1$
and $3-2\cos(l\varphi_0)\geq1,$ the equation holds if and only if $\cos(k\varphi_0)=1$
and $\cos(l\varphi_0)=1.$ So $k\varphi_0=2\pi s_0$ and $\cos(l\varphi_0)=2\pi t_0$ for some
integer $s_0$ and $t_0.$  As $k$ and $l$ are relatively prime,  there exist two integers
$p$ and $q$ such that $kp+ql=1.$ Hence $\varphi_0=\varphi_0(kp+lq)=2\pi(ps_0+qt_0)\in2\pi\mathbb{Z}.$
As a result, $z_0=e^{i\,\varphi_0}=1.$ Now we have to show that the multiplicity of the root
$z_0=1$ is $2.$ Indeed, $P(1)=P^\prime(1)=0$ and $P^{\prime\prime}(1)=-2(k^2+l^2)\neq0.$
$\hfill \qed$

Let us set $H(z)=\prod\limits_{s=1}^{m}(z-z_s)(z-z_s^{-1}),$ where $m=k+l-1$ and $z_s$ are roots
of $P(z)$ different from $1.$ Then by Lemma \ref{lemma2}, we have $P(z)=\frac{(z-1)^2}{z^{k+l}} H(z).$

\medskip

\begin{lemma}\label{lemma5}  Let $H(z)=\prod\limits_{s=1}^{m}(z-z_s)(z-z_s^{-1})$ and  $H(1)\neq 0.$ Then
$$\prod\limits_{j=1}^{n-1}H(\varepsilon_n^j)=\prod\limits_{s=1}^{m}\frac{T_n(w_s)-1}{w_s-1},$$  where
$w_s=\frac12(z_s+z_s^{-1}),\,s=1,\ldots,m$  and $T_n(x)$  is the
Chebyshev polynomial of the first kind.
\end{lemma}

\textbf{Proof.} It is easy to check that $\prod\limits_{j=1}^{n-1}(z-\varepsilon_n^j)=\frac{z^n-1}{z-1}$
if $z\neq1.$ Also we note that $\frac12(z^n+z^{-n})=T_n(\frac12(z+z^{-1})).$  By the substitution
$z=e^{i\,\varphi},$ the latter follows from the evident identity $\cos(n\varphi)=T_n(\cos\varphi).$
Then we have
\begin{eqnarray*}
\prod\limits_{j=1}^{n-1} H(\varepsilon_n^j) &=&
\prod\limits_{j=1}^{n-1}\prod\limits_{s=1}^{m}(\varepsilon_n^j-z_s)(\varepsilon_n^j-z_s^{-1})\\
&=&\prod\limits_{s=1}^{m}\prod\limits_{j=1}^{n-1}(z_s-\varepsilon_n^j)(z_s^{-1}-\varepsilon_n^j)\\
&=&\prod\limits_{s=1}^{m}\frac{z_s^n-1}{z_s-1}\frac{z_s^{-n}-1}{z_s^{-1}-1}=
\prod\limits_{s=1}^{m}\frac{T_n(w_s)-1}{w_s-1}.
\end{eqnarray*}
$\hfill \qed$

Note that
$\prod\limits_{j=1}^{n-1}(1-\varepsilon_n^j)=\lim\limits_{z\to1}\prod\limits_{j=1}^{n-1}
(z-\varepsilon_n^j)= \lim\limits_{z\to1}\frac{z^n-1}{z-1}=n$ and
$\prod\limits_{j=1}^{n-1}\varepsilon_n^{j} = (-1)^{n-1}$.  As a
result, taking into account Lemma~\ref{lemma2} and Lemma~\ref{lemma5}, we obtain
\begin{eqnarray*}
\tau_{k,l}(n)&=&\frac{1}{n}\prod\limits_{j=1}^{n-1}P(\varepsilon_n^j)=\frac{1}{n}\prod\limits_{j=1}^{n-1}
\frac{(\varepsilon_n^j-1)^2}{(\varepsilon_n^{j})^{k+l}}H(\varepsilon_n^j)=\frac{(-1)^{(n-1)(k+l)}n^2}{n}
\prod\limits_{j=1}^{n-1}H(\varepsilon_n^j)\\
&=&(-1)^{(n-1)(k+l)}n\prod\limits_{s=1}^{k+l-1}\frac{T_n(w_s)-1}{w_s-1}.
\end{eqnarray*}
$\hfill \qed$

\begin{cor}\label{corollary2}
$\tau_{k,l}(n)=n \left|\prod_{s=1}^{k+l-1}U_{n-1}(\sqrt{\frac{1+w_s}{2}})\right|^2,$
where $w_s, s=1,2,\ldots,k$  are the same as in Theorem~\ref{maintheorem}
and $U_{n-1}(w)$ is the Chebyshev polynomial of the second kind.
\end{cor}

\smallskip
\textbf{Proof.} Follows from the identity
$\frac{T_n(w)-1}{w-1}=U_{n-1}^2(\sqrt{\frac{1+w }{2}}).  \hfill \qed $

The following theorem appeared after a fruitful discussion with professor D. Lorenzini.

\begin{theorem}\label{Theortausquare}

Let $\tau(n)=\tau_{k,l}(n)$ be the number of spanning trees of the graph $I(n,k,l).$
Then there exist an integer sequence $a(n)=a_{k,l}(n), n\in\mathbb{N}$ such that
\begin{itemize}
\item[$1^\circ$] $\tau(n)=n\,a^2(n)$ when $n$ is odd
\item[$2^\circ$] $\tau(n)=6n\,a^2(n)$ when $n$ is even and $k+l$ is even,
\item[$3^\circ$] $\tau(n)=n\,a^2(n)$ when $n$ is even and $k+l$ is odd.
\end{itemize}
\end{theorem}

\textbf{Proof.}
Recall that all nonzero eigenvalues of $I(n,k,l)$ are given by the list
$\{\lambda_{2,0},\,\lambda_{1,j},\lambda_{2,j},\,j=1,\ldots, n-1\}.$
By the Kirchhoff theorem we have
$2n\tau(n)=\lambda_{2,0}\prod_{j=1}^{n-1}\lambda_{1,j}\lambda_{2,j}.$

Since $\lambda_{2,0}=2,$ we have $n\tau(n)=\prod_{j=1}^{n-1}\lambda_{1,j}\lambda_{2,j}.$
We note that
$\lambda_{1,j}\lambda_{2,j}=P(\varepsilon_n^j)=P(\varepsilon_n^{n-j})=\lambda_{1,n-j}\lambda_{2,n-j}.$
So, we get $n\tau(n)=(\prod_{j=1}^{(n-1)/2}\lambda_{1,j}\lambda_{2,j})^2$ if $n$ is odd and
$n\tau(n)=\lambda_{1,\frac{n}{2}}\lambda_{2,\frac{n}{2}}(\prod_{j=1}^{n/2-1}\lambda_{1,j}\lambda_{2,j})^2,$
if $n$ is even. The value $\lambda_{1,\frac{n}{2}}\lambda_{2,\frac{n}{2}}=P(-1)=(3-2(-1)^{k})(3-2(-1)^l)-1$
is equal to $24,$ if $k$ and $l$ are of different parity and $4,$ if both $k$ and $l$ are odd.
The case when both $k$ and $l$ are even is impossible, since $k$ and $l$ are relatively prime.

The graph $I(n,k,l)$ admits a cyclic group of automorphisms isomorphic to $\mathbb{Z}_n,$ which
acts freely on the set of spanning trees. Therefore, the value $\tau(n)$ is a multiple of $n$.
So $\frac{\tau(n)}{n}$ is an integer. Hence
\begin{itemize}
\item[$1^\circ$] $\frac{\tau(n)}{n}=(\frac{\prod_{j=1}^{(n-1)/2}\lambda_{1,j}\lambda_{2,j}}{n})^2$
when $n$ is odd,
\item[$2^\circ$] $\frac{\tau(n)}{n}=(\frac{2\prod_{j=1}^{n/2-1}\lambda_{1,j}\lambda_{2,j}}{n})^2$
when $n$ is even and $k+l$ is even,
\item[$3^\circ$] $\frac{\tau(n)}{n}=6(\frac{2\prod_{j=1}^{n/2-1}\lambda_{1,j}\lambda_{2,j}}{n})^2$
when $n$ is even and $k+l$ is odd.
\end{itemize}

Each algebraic number $\lambda_{i,j}$ comes into both products
$\prod_{j=1}^{(n-1)/2}\lambda_{1,j}\lambda_{2,j}$ and $\prod_{j=1}^{n/2-1}\lambda_{1,j}\lambda_{2,j}$
with all its Galois conjugate elements. Therefore, both products are integer numbers. From here we
conclude that in equalities $1^\circ,\,2^\circ$ and $3^\circ$ the value that is squared is a rational
number. Because $\frac{\tau(n)}{n}$ is integer and $6$ is a squarefree, all these rational numbers
are integer. Setting $a(n)=\frac{\prod_{j=1}^{(n-1)/2}\lambda_{1,j}\lambda_{2,j}}{n}$ if $n$ is odd
and $a(n)=\frac{2\prod_{j=1}^{n/2-1}\lambda_{1,j}\lambda_{2,j}}{n}$ if $n$ is even, we finish the proof
of the theorem.
\qed

\bigskip

From now on, we aim to  estimate the minimum number of generators
for the Jacobian of $I$-graph $I(n,k,l).$

\begin{lemma}\label{spantreebign3}
For any given $I$-graph $I(n,k,l)$ the number of spanning trees
$\tau(n)$ satisfies the inequality $\tau(n)\geq n^3.$
\end{lemma}

\begin{pf}
Recall that for any $j=0,\ldots, n-1,$ the Laplacian matrix $L$ of
$I(n,k,l)$ has two eigenvalues, say $\lambda_{1,j}$ and
$\lambda_{2,j},$ which are roots of the quadratic equation
$Q_j(\lambda) =(3-\lambda-\varepsilon_n^{jk}-\varepsilon_n^{-jk})
(3-\lambda-\varepsilon_n^{jl}-\varepsilon_n^{-jl})-1=0.$ So,
$\lambda_{1,
j}\lambda_{2,j}=(3-\varepsilon_n^{jk}-\varepsilon_n^{-jk})
(3-\varepsilon_n^{jl}-\varepsilon_n^{-jl})-1=P(\varepsilon_n^j).$
Note that $\lambda_{1,0}=0$ and $\lambda_{2,0}=2.$ Furthermore,
$\{\lambda_{1,j},\lambda_{2,j}\,|\,j=0,\ldots, n-1\}$ is the set of
all eigenvalues of $L.$ The Kirchhoff theorem states the following
$$2n\,\tau_{k,l}(n)=2n\,\tau(n)=\lambda_{2,0}\prod_{j=1}^{n-1}\lambda_{1,j}\lambda_{2,j}
=2\prod_{j=1}^{n-1}\lambda_{1,j}\lambda_{2,j}.$$

Hence $n\tau(n)=\prod_{j=1}^{n-1}P(\varepsilon_n^j),$ where
$P(\varepsilon_n^j)=(3-2\cos(\frac{2 j k\pi}{n}))(3-2\cos(\frac{2 j
l\pi}{n}))-1.$ It is easy to prove the following trigonometric
identity
$$(3-2\cos(\frac{2 j k\pi}{n}))(3-2\cos(\frac{2 j l\pi}{n}))-1=
4\sin(\frac{j k\pi}{n})^2+4\sin(\frac{j l\pi}{n})^2+16\sin(\frac{j
k\pi}{n})^2\sin(\frac{j l\pi}{n})^2.$$

Connectedness of $I$-graph implies $GCD(n,k,l)=1.$ It may happen
that $GCD(n,k)=m\neq1$ and $GCD(n,l)=m'\neq1.$ We will use the
notation $n=m\,q=m'q',\,k=p\,m,\,l=p'm'.$ We introduce three sets,
$J,\,J_k$ and $J_l$ in the following way
$J=\{1,2,\ldots,n-1\},\,J_k=\{j\big|j=d\,q,\,d=1,\ldots,m-1\}$ and
$J_l=\{j\big|j=d\,q',\,d=1,\ldots,m'-1\}.$ If $j\in J_k$ then
$\sin(\frac{j\,k\,\pi}{n})=0$ and if $j\in J_l$ then
$\sin(\frac{j\,l\,\pi}{n})=0.$ Now we are going to find a low bound
for $\tau(n).$ As $n\,\tau(n)=\prod_{j=1}^{n-1}P(\varepsilon_n^j),$
we evaluate the product.

\begin{eqnarray*}
\prod_{j=1}^{n-1}P(\varepsilon_n^j)&=&\prod_{j=1}^{n-1}(4\sin(\frac{j
k\pi}{n})^2+
4\sin(\frac{j l\pi}{n})^2+16\sin(\frac{j k\pi}{n})^2\sin(\frac{j l\pi}{n})^2)\\
&\geq&\prod_{j\in J_k}4\sin(\frac{j l\pi}{n})^2\prod_{j\in
J_l}4\sin(\frac{j k\pi}{n})^2
\prod_{j\in J\setminus(J_k\cup J_l)}16\sin(\frac{j k\pi}{n})^2\sin(\frac{j l\pi}{n})^2\\
&=&\prod_{j\in J\setminus J_k}4\sin(\frac{j k\pi}{n})^2 \prod_{j\in
J\setminus J_l}4\sin(\frac{j l\pi}{n})^2.
\end{eqnarray*}

Now we analyze individual components of the product. We make use of
the following simple identity $\cos(\frac{2 j
p\pi}{q})=\cos(\frac{2(j+q)p\pi}{q}).$

\begin{eqnarray*}
\prod_{j\in J\setminus J_k}4\sin(\frac{j k\pi}{n})^2&=& \prod_{j\in
J\setminus J_k}(2-2\cos(\frac{2 j k\pi}{n}))=
\prod_{j\in J\setminus J_k}(2-2\cos(\frac{2 j m p\pi}{m q}))\\
&=&\prod_{j\in J\setminus J_k}(2-2\cos(\frac{2 j p\pi}{q}))=
\prod_{j=1}^{q-1}(2-2\cos(\frac{2 j p\pi}{q}))^m.
\end{eqnarray*}

The Chebyshev polynomial $T_q(x)=\cos(q\arccos(x))$ has the
following property. The roots of the equation $T_q(x)-1=0$ are
$\cos(\frac{2 j \pi}{q}),\,j=0,1,\ldots,q-1.$ Since the leading
coefficient of $T_q(x)$ is $2^{q-1},$ for $x\neq1$ we have the
identity

$$\prod_{j=1}^{q-1}(2x-2\cos(\frac{2 j \pi}{q}))=\frac{T_q(x)-1}{x-1}.$$

As $p$ and $q$ are co-prime, we obtain

$$\prod_{j=1}^{q-1}(2-2\cos(\frac{2 j p\pi}{q}))^m=\prod_{j=1}^{q-1}(2-2\cos(\frac{2 j \pi}{q}))^m
=(\lim_{x\rightarrow1}\frac{T_q(x)-1}{x-1})^m=(q^2)^m=(\frac{n}{m})^{2m}.$$

Hence

$$\prod_{j\in J\setminus J_k}4\sin(\frac{j k\pi}{n})^2=\big(\frac{n}{m}\big)^{2m}.$$

In a similar way, we obtain
$$\prod_{j\in J\setminus J_l}4\sin(\frac{j l\pi}{n})^2=\big(\frac{n}{m'}\big)^{2m'}.$$

To get the final result we use the following trivial inequality. For
any integers $a\geq2$ and $b\geq2$ we have $a^b\geq ab.$ Since
$q=n/m\geq2$ and $q'=n/m'\geq2,$ we conclude
$$n\,\tau(n)=\prod_{j=1}^{n-1}P(\varepsilon_n^j)\geq(\frac{n}{m})^{2m}(\frac{n}{m'})^{2m'}\geq n^2 n^2=n^4.$$
\end{pf}

Using Lemma \ref{spantreebign3}, one can show the following theorem.

\begin{theorem}\label{GenRankIGr}
For any given $I$-graph $I(n,k,l)$ the minimum number of generators
for Jacobian $Jac(I(n,k,l))$ is at least $2$ and at most $2k+2l-1.$
The both bounds are sharp.
\end{theorem}

\begin{pf}
 The upper bound for
the number of generators follows from theorem~\ref{theorem0}.
Indeed, by this theorem the group $coker(L(I(n,k,l))\cong
Jac(I(n,k,l))\oplus\mathbb{Z}$ is generated by $2k+2l$ elements. One
of these generators is needed to generate the infinite cyclic group
$\mathbb{Z}.$ Hence $Jac(I(n,k,l))$ is generated by $2k+2l-1$
elements.

To get the lower bound we use Lemma~\ref{spantreebign3}. Let us
suppose that $Jac(I(n,k,l))$ is generated by one element. Then it is
the cyclic group of order $\tau(n).$ Denote by $D$ be a product of
all distinct nonzero eigenvalues of $I(n,k,l)$ By Proposition 2.6
from \cite{Lor}, the order of each element of $Jac(I(n,k,l))$ is
divisor of $D.$ Hence, $\tau(n)$ is divisor of $D$ and we have
inequality $D\geq\tau(n).$ By the Kirchhoff theorem we have
$2n\tau(n)=\lambda_{2,0}\prod_{j=1}^{n-1}\lambda_{1,j}\lambda_{2,j}.$
We note that any algebraic number $\lambda_{i,j}$ comes into the
product together with its Galois conjugate, so $2n\tau(n)$ is a
multiple of $D.$ In particular, $2n\tau(n)\geq D.$

From the proof of Theorem~\ref{Theortausquare} we have $n\tau(n)=(\prod_{j=1}^{(n-1)/2}\lambda_{1,j}\lambda_{2,j})^2,$
if $n$ is odd and $n\tau(n)=\lambda_{1,\frac{n}{2}}\lambda_{2,\frac{n}{2}}
(\prod_{j=1}^{n/2-1}\lambda_{1,j}\lambda_{2,j})^2,$ if $n$ is even. Moreover, the value
$\lambda_{1,\frac{n}{2}}\lambda_{2,\frac{n}{2}}$ is equal to $24,$ if $k$ and $l$ are of
different parity and $4,$ if both $k$ and $l$ are odd. The case when both $k$ and $l$ are
even is impossible as $k$ and $l$ are relatively prime.

Now, we have $4n\tau(n)=(2\prod_{j=1}^{(n-1)/2}\lambda_{1,j}\lambda_{2,j})^2$ if $n$ is
odd. We note that any algebraic number $\lambda_{i,j}$ comes into the product
$\rho=2\prod_{j=1}^{(n-1)/2}\lambda_{1,j}\lambda_{2,j}$ together with its Galois conjugate.
Therefore, the product $\rho$ is an integer number and contains all distinct nonzero eigenvalues.
Hence $\rho$ is a multiple of $D.$ So, we obtain $4n\tau(n)=\rho^2\geq D^2\geq\tau(n)^2.$

Also we get
$4n\lambda_{1,\frac{n}{2}}\lambda_{2,\frac{n}{2}}\tau(n)=
(2\lambda_{1,\frac{n}{2}}\lambda_{2,\frac{n}{2}}\prod_{j=1}^{n/2-1}\lambda_{1,j}\lambda_{2,j})^2$
if $n$ is even. By a similar argument, taking into account the
inequality $24\geq\lambda_{1,\frac{n}{2}}\lambda_{2,\frac{n}{2}}$ we
obtain
$96n\tau(n)\geq4n\lambda_{1,\frac{n}{2}}\lambda_{2,\frac{n}{2}}\tau(n)\geq
D^2\geq\tau(n)^2.$

As result, by Lemma~\ref{spantreebign3} we have $4n\geq\tau(n)\geq
n^3$ if $n$ is odd and $96n\geq\tau(n)\geq n^3$ if $n$ is even. For
$n\geq10$ this is impossible. So, the rank of $Jac(I(n,k,l))$ is at
least two for all $n\geq10.$ For $n$ less than $10$ this statement
can be proved by direct calculation.
\end{pf}
The upper bound $2k+2l-1$ for the number of generators of
$Jac(I(n,k,l))$ is attained for graphs $I(17,2,3)$ and $I(170,3,4).$
See Tables $2$ and $3$ in section 7.

\section{Asymptotic for the number of spanning trees}

The asymptotic for the number of spanning trees of the graph $I(n,k,l)$ is given in the following theorem.

\begin{theorem}
Let $P(z)=(3-z^{k}-z^{-k})(3-z^{l}-z^{-l})-1.$ Suppose that $k$ and $l$ are relatively prime and
set $A_{k,l}=\prod\limits_{P(z)=0,\,|z|>1}|z|.$ Then the number $\tau_{k,l}(n)$ of spanning trees
of the graph $I(n,k,l)$ has the asymptotic
$$\tau_{k,l}(n)\sim \frac{n}{k^2+l^2}A_{k,l}^n,\,n\to\infty.$$\end{theorem}

\textbf{Proof.} By theorem \ref{maintheorem} we have
$$\tau_{k,l}(n)=(-1)^{(n-1)(k+l)}n\prod_{s=1}^{k+l-1}\frac{T_n(w_s)-1}{w_s-1},$$ where
$w_s, s=1,2,\ldots,k+l-1$ are roots of the polynomial $Q(w)=\frac{(3-2T_k(w))(3-2T_l(w))-1}{w-1}.$
So we obtain
$$\tau_{k,l}(n)=n\prod_{s=1}^{k+l-1}\left|\frac{T_n(w_s)-1}{w_s-1}\right|=
n\prod_{s=1}^{k+l-1}|T_n(w_s)-1|\big/\prod_{s=1}^{k+l-1}|w_s-1|.$$

By lemma~\ref{lemma2} we have $T_n(w_s)=\frac{1}{2}(z_s^n+z_s^{-n}),$ where the $z_s$ and $1/z_s$ are roots
of the polynomial $P(z)$ with the property $|z_s|\neq1,\,s=1,2,\ldots,k+l-1.$ Replacing $z_s$ by $1/z_s,$ if
it is necessary, we can assume that all $|z_s|>1$ for all $s=1,2,\ldots,k+l-1.$ Then
$T_n(w_s)\sim\frac{1}{2}z_s^n$ as $n$ tends to $\infty.$ So $|T_n(w_s)-1|\sim\frac{1}{2}|z_s|^n$ as $n\to\infty.$
Hence
$$\prod_{s=1}^{k+l-1}|T_n(w_s)-1|\sim\frac{1}{2^{k+l-1}}\prod_{s=1}^{k+l-1}|z_s|^n=
\frac{1}{2^{k+l-1}}\prod\limits_{P(z)=0,\,|z|>1}|z|^n=\frac{1}{2^{k+l-1}}A_{k,l}^n.$$

Now we directly evaluate the quantity $\prod_{s=1}^{k+l-1}|w_s-1|.$ We note that
$Q(w)=a_0w^{k+l-1}+a_12^{k+l-2}+\ldots+a_{k+l-2}w+a_{k+l-1}$ is an integer polynomial
with the leading coefficient $a_0=2^{k+l}.$ From here we obtain
$$\prod_{s=1}^{k+l-1}|w_s-1|=\prod_{s=1}^{k+l-1}|1-w_s|=|\frac{1}{a_0}Q(1)|
=\frac{2(k^2+l^2)}{2^{k+l}}=\frac{k^2+l^2}{2^{k+l-1}}.$$

Indeed,
$Q(1)=\lim\limits_{w\to1}\frac{(3-2T_k(w))(3-2T_l(w))-1}{w-1}=-2T_k^\prime(1)(3-2T_l(1))-2T_l^\prime(1)(3-2T_k(1))
=-2kU_{k-1}(1)(3-2T_l(1))-2lU_{l-1}(1)(3-2T_k(1))=-2(k^2+l^2)$ and $a_0=2^{k+l}.$

In order to get the statement of the theorem, we combine the above mentioned results. Then
$$\tau_{k,l}(n)\sim n\,\frac{A_{k,l}^n}{2^{k+l-1}}\Big/\frac{k^2+l^2}{2^{k+l-1}}=
\frac{n}{k^2+l^2}A_{k,l}^n\textrm{ as }n\to\infty.$$

\medskip
\textbf{Remark:} It was noted by professor A.~Yu.~Vesnin that constant $A_{k,l}$ coincides
with the Mahler measure of Laurent polynomial $P(z)=(3-z^{k}-z^{-k})(3-z^{l}-z^{-l})-1.$ It
gives a simple way to calculate $A_{k,l}$ using the following formula
$$A_{k,l}=e^{\int\limits_{0}^{1}\log|P(e^{2 \pi i t})|\textrm{d}t}.$$
See, for example, (\cite{EverWard}, p. 6) for the proof.

\medskip
The numerical values for $A_{k,l},$ where $k$ and $l$ are relatively prime numbers $1\leq k\leq l\leq9$
will be given in Table $1$ in the section 7.

\section{Examples and Tables}

\subsection{Examples}

\begin{enumerate}
\item[$1^\circ$] The Prism graph $I(n,1,1).$
We have the following asymptotic $\tau_{1,1}(n)=n(T_n(2)-1)\sim \frac{n}{2}(2+\sqrt{3})^n,\,n\to\infty.$
\item[$2^\circ$] The generalized Petersen graph $GP(n,2)=I(n,1,2).$
The the number of spanning trees (see \cite{KwonMedMed}) behaves like
$\tau_{1,2}(n)\sim \frac{n}{5}A_{1,2}^n,\,n\to\infty,$ where
$A_{1,2}=\frac{7+\sqrt{5}+\sqrt{38+14\sqrt{5}}}{4}\cong4.39026.$
\item[$3^\circ$] The smallest proper $I$-graph $I(n,2,3)$ has the following asymptotic for the number
of spanning trees $\tau_{2,3}(n)\sim\frac{n}{13}A_{2,3}^n,\,n\to\infty.$ Here $A_{2,3}\cong4.84199$
is a suitable root of the algebraic equation
$1-7x+13x^2-35x^3+161x^4-287x^5+241x^6-371x^7+577x^8-371x^9+241x^{10}-287x^{11}+161x^{12}-35x^{13}+13x^{14}-7x^{15}+x^{16}=0.$
\end{enumerate}

Here is the table for asymptotic constants $A_{k,l}$ for relatively prime numbers $1\leq k\leq l\leq9.$

\begin{table}[h]
\caption{Asymptotic constants $A_{k,l}$}
\begin{tabular}{|c|c|c|c|c|c|c|c|c|c|}\hline
$k\backslash l$ & $1$ &$2$ &$3$ &$4$ &$5$ &$6$ &$7$& $8$& $9$\\ \hline
$1$ & $3.7320$ & $4.3902$ & $4.7201$ & $4.8954$ & $4.9953$ & $5.0559$ & $5.0945$ & $5.1203$ & $5.1382$\\ \hline
$2$ & & - & $4.8419$ & - & $5.0249$ & - & $5.1033$ & - & $5.1414$\\ \hline
$3$ & & & - & $5.0054$ & $5.0541$ & - & $5.1137$ & $5.1320$ & -\\ \hline
$4$ & & & & - & $5.0802$ & - & $5.1244$ & - & $5.1504$ \\ \hline
$5$ & & & & & - & $5.1201$ & $5.1346$ & $5.1461$ & $5.1554$\\ \hline
$6$ & & & & & & - & $5.1438$ & - & -\\ \hline
$7$ & & & & & &  & - & $5.1589$ & $5.1649$ \\ \hline
$8$ & & & & & &  &  & - & $5.1691$ \\ \hline
\end{tabular}
\end{table}

\subsection{The tables of Jacobians of $I$-graphs}

Theorem~\ref{theorem0} is the first step to understand the structure of the Jacobian for $I(n,k,l).$
Also, it gives a simple way for numerical calculations of $\textrm{Jac}(I(n,k,l))$ for small values
of $k$ and $l.$ See Tables $2$ and $3$ below.

\begin{table}[h]
\caption{Graph $I(n,2,3)$}
\begin{tabular}{r|l|r}
$n$ & $\textrm{Jac}(I(n,2,3))$ &$\tau_{2,3}(n) =|\textrm{Jac}(I(n,2,3))|$\\ \hline
4&$\mathbb{Z}_{7}\oplus\mathbb{Z}_{28}$& 196 \\
5&$\mathbb{Z}_{19}\oplus\mathbb{Z}_{95}$& 1805\\
6&$\mathbb{Z}_{19}\oplus\mathbb{Z}_{114}$& 2166\\
7&$\mathbb{Z}_{83}\oplus\mathbb{Z}_{581}$& 48223\\
8&$\mathbb{Z}_{161}\oplus\mathbb{Z}_{1288}$& 207368\\
9&$\mathbb{Z}_{289}\oplus\mathbb{Z}_{2601}$& 751689\\
10&$\mathbb{Z}_{1558}\oplus\mathbb{Z}_{3895}$& 6068410\\
11&$\mathbb{Z}_{1693}\oplus\mathbb{Z}_{18623}$& 31528739\\
12&$\mathbb{Z}_{5}\oplus\mathbb{Z}_{5}\oplus\mathbb{Z}_{665}\oplus\mathbb{Z}_{7980}$& 132667500\\
13&$\mathbb{Z}_{25}\oplus\mathbb{Z}_{325}\oplus\mathbb{Z}_{325}\oplus\mathbb{Z}_{325}$& 858203125\\
14&$\mathbb{Z}_{17513}\oplus\mathbb{Z}_{245182}$& 4293872366\\
15&$\mathbb{Z}_{37069}\oplus\mathbb{Z}_{556035}$& 20611661415\\
16&$\mathbb{Z}_{84847}\oplus\mathbb{Z}_{1357552}$& 115184214544\\
17&$\mathbb{Z}_{2}\oplus\mathbb{Z}_{2}\oplus\mathbb{Z}_{2}\oplus\mathbb{Z}_{2}
\oplus\mathbb{Z}_{2}\oplus\mathbb{Z}_{2}\oplus\mathbb{Z}_{23186}\oplus\mathbb{Z}_{394162}$& 584898568448\\
18&$\mathbb{Z}_{400843}\oplus\mathbb{Z}_{7215174}$& 2892151991682\\
19&$\mathbb{Z}_{898243}\oplus\mathbb{Z}_{17066617}$& 15329969253931\\
20&$\mathbb{Z}_{19}\oplus\mathbb{Z}_{19}\oplus\mathbb{Z}_{19}\oplus\mathbb{Z}_{19}\oplus\mathbb{Z}_{5453}
\oplus\mathbb{Z}_{109060}$& 77502443441780\\
21&$\mathbb{Z}_{4301807}\oplus\mathbb{Z}_{90337947}$& 388616412770229\\
22&$\mathbb{Z}_{9536669}\oplus\mathbb{Z}_{209806718}$& 2000857223542342\\
23&$\mathbb{Z}_{20949827}\oplus\mathbb{Z}_{481846021}$& 10094590780588367\\
24&$\mathbb{Z}_{5}\oplus\mathbb{Z}_{5}\oplus\mathbb{Z}_{9192295}\oplus\mathbb{Z}_{220615080}$& 50598972420215000\\
25&$\mathbb{Z}_{101468531}\oplus\mathbb{Z}_{2536713275}$& 257396569582449025\\
26&$\mathbb{Z}_{25}\oplus\mathbb{Z}_{325}\oplus\mathbb{Z}_{8923525}\oplus\mathbb{Z}_{17847050}$& 1293976099416406250\\
27&$\mathbb{Z}_{490309597}\oplus\mathbb{Z}_{13238359119}$& 6490894524578165043\\
28&$\mathbb{Z}_{49}\oplus\mathbb{Z}_{154342069}\oplus\mathbb{Z}_{4321577932}$& 32683062689111444092\\
29&$\mathbb{Z}_{2376466133}\oplus\mathbb{Z}_{68917517857}$& 163780147157583236981\\
30&$\mathbb{Z}_{19}\oplus\mathbb{Z}_{19}\oplus\mathbb{Z}_{275089049}\oplus\mathbb{Z}_{8252671470}$&
819549256247415262830\\
31&$\mathbb{Z}_{11507960491}\oplus\mathbb{Z}_{356746775221}$& 4105427794534925793511\\
32&$\mathbb{Z}_{25318259953}\oplus\mathbb{Z}_{810184318496}$& 20512457185525873990688\\
33&$\mathbb{Z}_{55700389051}\oplus\mathbb{Z}_{1838112838683}$& 102383600234281102459833\\
34&$\mathbb{Z}_{2}\oplus\mathbb{Z}_{4}\oplus\mathbb{Z}_{4}\oplus\mathbb{Z}_{4}\oplus\mathbb{Z}_{4}
\oplus\mathbb{Z}_{4}\oplus\mathbb{Z}_{4}\oplus\mathbb{Z}_{1915580948}\oplus\mathbb{Z}_{32564876116}$&
511022336096582352633856\\
35&$\mathbb{Z}_{269747901677}\oplus\mathbb{Z}_{9441176558695}$& 2546737566070056079431515\\
\end{tabular}
\end{table}

\clearpage
\begin{table}[h]
\caption{Graph $I(n,3,4)$}
\begin{tabular}{r|l|r} $n$ & $\textrm{Jac}(I(n,3,4))$ & $\tau_{3,4}(n)=|\textrm{Jac}(I(n,3,4))|$ \\ \hline
5& $\mathbb{Z}_{2}\oplus\mathbb{Z}_{10}\oplus\mathbb{Z}_{10}\oplus\mathbb{Z}_{10}$ & 2000 \\
6& $\mathbb{Z}_{19}\oplus\mathbb{Z}_{114}$ & 2166 \\
7& $\mathbb{Z}_{71}\oplus\mathbb{Z}_{497}$ & 35287 \\
8& $\mathbb{Z}_{73}\oplus\mathbb{Z}_{584}$ & 42632 \\
9& $\mathbb{Z}_{289}\oplus\mathbb{Z}_{2601}$ & 751689 \\
10& $\mathbb{Z}_{2}\oplus\mathbb{Z}_{12}\oplus\mathbb{Z}_{60}\oplus\mathbb{Z}_{60}\oplus\mathbb{Z}_{60}$ & 5184000 \\
11& $\mathbb{Z}_{1541}\oplus\mathbb{Z}_{16951}$ & 26121491 \\
12& $\mathbb{Z}_{11}\oplus\mathbb{Z}_{11}\oplus\mathbb{Z}_{209}\oplus\mathbb{Z}_{2508}$ & 63424812 \\
13& $\mathbb{Z}_{5}\oplus\mathbb{Z}_{5}\oplus\mathbb{Z}_{1555}\oplus\mathbb{Z}_{20215}$ & 785858125 \\
14& $\mathbb{Z}_{16969}\oplus\mathbb{Z}_{237566}$ & 4031257454 \\
15& $\mathbb{Z}_{2}\oplus\mathbb{Z}_{10}\oplus\mathbb{Z}_{17410}\oplus\mathbb{Z}_{52230}$ & 18186486000 \\
16& $\mathbb{Z}_{71321}\oplus\mathbb{Z}_{1141136}$ & 81386960656 \\
17& $\mathbb{Z}_{2}\oplus\mathbb{Z}_{2}\oplus\mathbb{Z}_{2}\oplus\mathbb{Z}_{2}\oplus\mathbb{Z}_{2}
\oplus\mathbb{Z}_{2}\oplus\mathbb{Z}_{23186}\oplus\mathbb{Z}_{394162}$ & 584898568448 \\
18& $\mathbb{Z}_{400843}\oplus\mathbb{Z}_{7215174}$ & 2892151991682 \\
19& $\mathbb{Z}_{37}\oplus\mathbb{Z}_{37}\oplus\mathbb{Z}_{23939}\oplus\mathbb{Z}_{454841}$ & 14906272578931 \\
20& $\mathbb{Z}_{8}\oplus\mathbb{Z}_{12}\oplus\mathbb{Z}_{120}\oplus\mathbb{Z}_{79080}\oplus\mathbb{Z}_{79080}$ & 72042006528000 \\
21& $\mathbb{Z}_{4487981}\oplus\mathbb{Z}_{94247601}$ & 422981442583581 \\
22& $\mathbb{Z}_{10002631}\oplus\mathbb{Z}_{220057882}$ & 2201157792287542 \\
23& $\mathbb{Z}_{22138559}\oplus\mathbb{Z}_{509186857}$ & 11272663275719063 \\
24& $\mathbb{Z}_{187}\oplus\mathbb{Z}_{187}\oplus\mathbb{Z}_{259369}\oplus\mathbb{Z}_{6224856}$ & 56458663080288216 \\
25& $\mathbb{Z}_{2114}\oplus\mathbb{Z}_{52850}\oplus\mathbb{Z}_{52850}\oplus\mathbb{Z}_{52850}$ & 312061332000250000 \\
\end{tabular} \vspace{1cm}
\end{table}

The first example of Jacobian $Jac(I(n,3,4))$ with the maximum rank 13:
$$n=170,$$
$$Jac(I(170,3,4))\cong\mathbb{Z}_{2}\oplus\mathbb{Z}_{4}^8\oplus\mathbb{Z}_{6108}\oplus\mathbb{Z}_{30540}
\oplus\mathbb{Z}_{2^{2}\cdot3\cdot5\cdot103\cdot509\cdot1699\cdot11593\cdot p\cdot q}
\oplus\mathbb{Z}_{2^{2}\cdot3\cdot5\cdot17\cdot103\cdot509\cdot1699\cdot11593\cdot p\cdot q},$$
and
$$\tau_{3,4}(170)=2^{25}\cdot3^{4}\cdot5^{3}\cdot17\cdot103^{2}\cdot509^{4}\cdot1699^{2}\cdot
11593^{2}\cdot p^{2}\cdot q^{2},$$
where $p=16901365279286026289$ and $q=34652587005966540929.$

\clearpage
\newpage
\section*{ACKNOWLEDGMENTS}

The author is grateful to professor D. Lorenzini for helpful comments on the preliminary results of the paper
and professor Young Soo Kwon, whose remarks and suggestions assisted greatly in completion of the text.

The author was supported by the Russian Foundation for Basic Research (16-31-00138)
and the Slovenian-Russian grant (2016-2017).


\begin{thebibliography}{99}

\bibitem{CoriRoss} R. Cori, D. Rossin,
On the sandpile group of dual graphs,
European J. Combin. 21(4), (2000), 447--459.

\bibitem{BakerNorine} B. Baker, S. Norine,
Harmonic morphisms and hyperelliptic graphs,
Int. Math. Res. Notes. 15, (2009), 2914--2955.

\bibitem{Biggs} N.L. Biggs,
Chip-firing and the critical group of a graph,
J. Algebraic Combin. 9(1), (1999), 25--45.

\bibitem{BachHarpNagnib} R. Bacher, P. de la Harpe, T. Nagnibeda,
The lattice of integral flows and the lattice of integral cuts on a finite graph,
Bull. Soc. Math. France. 125, (1997), 167--198.

\bibitem{BoePro} F.T. Boesch, H. Prodinger,
Spanning tree formulas and Chebyshev polynomials,
Graphs and Combinatorics. 2(1), (1986), 191--200.

\bibitem{GS11} R. Gera and P. St$\breve{\textrm{a}}$nic$\breve{\textrm{a}}$,
The spectrum of generalized Petersen graphs,
Australas. J. Combin. 49, (2011), 39--45.

\bibitem{FGW71} R. Frucht, J.E. Graver and M.E. Watkins,
The groups of the generalized Petersen graphs,
Proc. Cambridge Philos. Soc. 70, (1971), 211--218.

\bibitem{Lor} D. Lorenzini,
Smith normal form and Laplacians,
J. Combin. Theory Ser. B. 98(6), (2008), 1271--1300.

\bibitem{ZhangYongGol} Zhang Yuanping, Yong Xuerong, M.J. Golin,
The number of spanning trees in circulant graphs,
Discrete. Math. 223(1), (2000), 337--350.

\bibitem{ZhangYongGolin} Zhang Yuanping, Xuerong Yong, M.J. Golin,
Chebyshev polynomials and spanning tree formulas for circulant and related graphs,
Discrete Math. 298(1), (2005), 334--364.

\bibitem{XiebinLinZhang} Chen Xiebin, Qiuying Lin, Fuji Zhang,
The number of spanning trees in odd valent circulant graphs,
Discrete Math. 282(1), (2004), 69--79.

\bibitem{NP04} S.D. Nikolopoulos and C. Papadopoulos,
The number of spanning trees in $K_n$-complements of quasi-threshold graphs,
Graph Combinator, 20, (2004), 383--397.

\bibitem{Dhar} D. Dhar, P. Ruelle, S. Sen, D.-N. Verma,
Algebraic aspects of abelian sandpile models,
J. Phys. A. 28, (1995), 805--831.

\bibitem{SW00} R. Shrock and F.Y. Wu,
Spanning trees on graphs and lattices in d-dimensions,
J. Phys. A 33, (2000), 3881--3902.

\bibitem{SWZ16} W. Sun, S. Wang and J. Zhang,
Counting spanning trees in prism and anti-prism graphs,
J. Appl. Anal. Comput. 6, (2016), 65--75.

\bibitem{KotaniSunada} M. Kotani, T. Sunada,
Jacobian tori associated with a finite graph and its abelian covering graphs,
Adv. Appl. Math. 24, (2000), 89--110.

\bibitem{YaoChinPen} Yaoping Hou, Chingwah Woo, Pingge Chen,
On the sandpile group of the square cycle,
Linear Algebra Appl. 418, (2006), 457--467.

\bibitem{ChenHou} Chen, Pingge and Hou, Yaoping,
On the critical group of the Mobius ladder graph,
Austral. J. Combin. 36, (2006), 133--142.

\bibitem{MedZind} I.A. Mednykh, M.A. Zindinova,
On the structure of picard group for moebius ladder,
Sib. Electron. Math. Rep. 8, (2011), 54--61.

\bibitem{MedMed2} A.D. Mednykh, I.A. Mednykh,
On the structure of the Jacobian group for circulant graphs,
Doklady Mathematics. 94(1), (2016), 445--449.

\bibitem{PJDav} P.J. Davis, Circulant Matrices, AMS Chelsea Publishing, 1994.

\bibitem{HorPisZit} B. Horvat, T. Pisanski, A. $\check{\textrm{Z}}$itnik,
Isomorphism checking of I-graphs,
Graphs and Combinatorics. 28(6), (2012), 823--830.

\bibitem{HorBou} J.D. Horton, I.Z. Bouwer,
Symmetric Y-graphs and H-graphs,
Journal of Combinatorial theory, Series B, 53, (1991), 114--129.

\bibitem{PetkZakr} M. Petkovs$\check{\textrm{e}}$k, H. Zakrajs$\check{\textrm{e}}$k,
Enumeration of I-graphs: Burnside does it again,
Ars Mathematica Contemporanea 2, (2009), 241--262.

\bibitem{BobPisZit} M. Boben, T. Pisanski, A. $\check{\textrm{Z}}$itnik,
I-graphs and the corresponding configurations,
Journal of Combinatorial Designs. 13(6), (2005), 406--424.

\bibitem{BoCheMoSta} I.Z. Bouwer, W.W. Chernoff, B. Monson, \& Star, Z. (1988).
The Foster Census. Charles Babbage Research Centre, Winnipeg.

\bibitem{KwonMedMed} Y.S. Kwon, A.D. Mednykh, I.A. Mednykh,
On Jacobian group and complexity of the generalized Petersen graph $GP(n, k)$ through Chebyshev polynomials,
arXiv preprint arXiv:1612.03372. – 2016.

\bibitem{OlivVina} A.S.S. de Oliveira, C. Vinagre,
The spectrum of an I-graph,
arXiv preprint arXiv:1511.03513. – 2015.

\bibitem{Biggs71} N.L. Biggs, Three remarkable graphs, Canad. J. Math. 25, (1973), 397--411.

\bibitem{FruGraWat} R. Frucht, J. E. Graver, and M. E. Watkins,
The groups of the generalized Petersen graphs,
Proc. Cambridge Philos. Sot. 70, (1971), 211--218.

\bibitem{EverWard} Everest G., Ward T.,
Heights of polynomials and entropy in algebraic dynamics.
Springer Science \& Business Media, 2013.

\end{thebibliography}
\end{document}